\documentclass[12pt]{article}
 
\usepackage[margin=1in]{geometry} 
\usepackage{setspace}
\usepackage{amsmath,amsthm,amssymb}
\usepackage{mathtools}
\usepackage{faktor}
\usepackage{xfrac}
\usepackage{dsfont}
\usepackage{bbm}
\usepackage{systeme}
\usepackage{mathtools}
\usepackage{tikz-cd}
\usetikzlibrary{cd}
\usepackage{hyperref}
\usepackage{cite}
\usepackage{blindtext}
\usepackage{url}
\usepackage{indentfirst}
\usepackage{mathrsfs}

\newcommand{\N}{\mathbb{N}}

\DeclareMathAlphabet{\mathbbe}{U}{bbold}{m}{n}

\newcommand{\Z}{\mathbb{Z}}

\newcommand{\defeq}{\vcentcolon=}

\newcommand{\spec}[1]{\text{Spec} #1}

\newcommand{\Fun}{\textit{Fun}}

\theoremstyle{plain}
\newtheorem{thm}{Theorem}[section]% reset theorem numbering for each section

\newtheorem{prop}[thm]{Proposition}

\theoremstyle{definition}
\newtheorem{defn}[thm]{Definition} % definition numbers are dependent on theorem numbers
\newtheorem{eg}[thm]{Example} % same for example numbers

\title{\textsc{A brief introduction to derived schemes}}
\author{Expository paper by Aurel Malapani\\
Johns Hopkins University\\
Supervised by Professor Emily Riehl}
\date{Spring 2018}

\begin{document}
\maketitle
\bibliographystyle{amsplain}

	\begin{abstract}
		The development of mathematics has been characterized by the increasing interconnectivity of seemingly separate disciplines. Such interplay has been facilitated by a massive development in formalism; category theory has provided a common language to the study of mathematics. In the opposite direction, we have seen algebraic topology and category theory synthesize via higher category theory. In this expository paper, we examine an application of higher category theory to geometry through the development of simplicial rings and derived schemes.
		
	\end{abstract}
	
	\begin{section}{Preliminaries}
	
	We briefly review the definition of  a scheme, as presented in \cite{gw}. We then go on to provide some intuition underlying the notions of moduli problems, moduli spaces, and stacks. We then review the basics of model category theory. The term ``ring'' will mean a commutative and unital ring. 
		\begin{subsection}{A bit of algebraic geometry}
		\begin{subsubsection}{Schemes}
		 For a ring $R,$ let $\spec R$ denote the prime spectrum of $R$ endowed with the \emph{Zariski topology} in which closed sets are of the form 
		\[V(\mathfrak a) \defeq \left\{\mathfrak p \in \spec A; \mathfrak a \subseteq \mathfrak p\right\}\] 
		for ideals $\mathfrak a \subseteq R.$ Open sets of the form 
		\[D(f) \defeq \left\{\mathfrak p \in \spec A; (f) \not\subseteq \mathfrak p: f \in R\right\}\]
		form a basis for the topology. Moreover, there is a natural presheaf (in fact a sheaf) $\mathscr{O}_{\spec R}$ of rings on open subsets of $\spec R$ given by mapping $D(f)$ to the localization of $R$ at $f.$ We obtain a ringed space $(\spec R, \mathscr{O}_R).$ For the sake of completeness, we recall the definition of a ringed space.
		\begin{defn}
			Consider an arbitrary topological space $X.$ Denote by $\mathbf{Ouv}_X$ the category whose objects are open subsets of $X$ and whose morphisms are the inclusions. A \emph{ringed space} is a pair $(X, \mathscr{O}_X)$ where $X$ is a topological space and $\mathscr{O}_X$ is a sheaf of rings on $\mathbf{Ouv}_X.$
			
			A morphism of ringed spaces consists of a pair $(f, f^\flat)\colon (X, \mathscr{O}_Y) \to (Y, \mathscr{O}_Y)$ such that $f\colon X \rightarrow Y$ is a continuous and $f^\flat \colon \mathscr{O}_Y \rightarrow \mathscr{O}_X$ is such that the restriction $f^\flat \colon \mathscr{O}_Y(U) \rightarrow \mathscr{O}_X\left(f^{-1}(U)\right)$ is a ring homomorphism for every open $U \subseteq Y.$
		\end{defn}
		
		 A ringed space is called an \emph{affine scheme} if it is isomorphic to $(\spec R, \mathscr{O}_R)$ for some commutative unital ring $R.$ 
		 
		 There is an equivalence of categories $\mathbf{Comm}^\text{op} \xrightarrow{\sim} \mathbf{Aff},$ sending $R \mapsto \spec R$ between the category of commutative unital rings and the category of affine schemes \cite[Theorem 2.35]{gw}. A \emph{scheme} is obtained by gluing together affine schemes. More precisely, a ringed space $(X, \mathscr{O}_X)$ is a scheme if it admits an open covering $X = \cup_i U_i$ by open affine subspaces. Morphisms of schemes are morphisms of ringed spaces. By the equivalence of categories $\mathbf{Comm}^\text{op} \xrightarrow{\sim} \mathbf{Aff},$ together with the fact that $\Z$ is initial in $\mathbf{Comm},$ we see that $\spec \Z$ is terminal in $\mathbf{Aff}.$ In fact, it is not difficult to see that $\spec \Z$ is the terminal object in $\mathbf{Sch},$ the category of schemes. By some mild abuse of notation, we may identify $\mathbf{Sch}$ with the slice category $\mathbf{Sch}/\spec \Z.$ In full generality, one often considers the slice category $\mathbf{Sch}/S$ of schemes over some fixed base-scheme $S.$ We can think of $X \xrightarrow{f} S \in \mathbf{Sch}/S,$ as a family of schemes (the fibers of $f$) parametrized by points in $S.$ Moreover, given some scheme $T,$ it is correct to view the set of morphisms $T \to X \in \mathbf{Sch}$ as the set of $T$-valued points on $X.$ More generally, we may consider some functor $F\colon \mathbf{Sch}^\text{op} \rightarrow \mathbf{Set}$ (not necessarily representable) as a ``geometric object'' and $F(T)$ as the $T$-valued points on $F.$ 
		
		Another notion that will be used below is that of a \emph{quasi-coherent} $\mathscr{O}_X$ module on a scheme $(X, \mathscr{O}_X).$ 
		\begin{defn}[{Quasi-coherent $\mathscr{O}_X$-modules \cite[Definitions 7.1 \& 7.15]{gw}}]
		A $\mathscr{O}_X$-module is a sheaf $\mathscr{F}$ on open subsets of $X$ so that for every open $U \subseteq X,$ we have that the section $\mathscr{F}(U)$ is an $\mathscr{O}_X(U)$-module. Moreover, $\mathscr{F}$ is \emph{quasi-coherent} if for every open affine $\spec R \subseteq X$ there is an $R$-module $M$ such that for every $f \in R$ there is an isomorphism $\mathscr{F}(D(f)) \cong M_f,$ where $M_f$ is the localization of $M$ at $f.$ 
		\end{defn}	
		
		\end{subsubsection}
		
		\begin{subsubsection}{Moduli Spaces and Stacks}
		\par 
		Our brief introduction to schemes above allows us to ask the following question: for an arbitrary base scheme $S,$ what do isomorphism classes of schemes parametrized by $S$ look like? This is a typical example of a \emph{moduli problem}. We give a very brief (and hence incomplete) description of moduli problems, moduli spaces, and stacks. 
		
		In the most basic case, a moduli problem (on $\mathbf{Sch}$) is a functor $F \colon \mathbf{Sch}^\text{op} \to \mathbf{Set}$ so that $F(S)$ consists of isomorphism classes of schemes parametrized by $S.$ A moduli space $M$ is a space that corepresents $F,$ i.e., $M$ is a (fine) moduli space for the functor $F$ if there is a natural isomorphism 
		\[\begin{tikzcd}[column sep = small] \alpha\colon F \ar[r,Rightarrow, "\cong"] & \text{Hom}(-, M)\end{tikzcd}.\] 
		In other words, borrowing from the concept of a covering space in algebraic topology, we require that $M$ has a universal family. In the language of higher category theory, this is also known as a ``universal element.''\cite{riehl}
		
		It is certainly not true that every moduli problem admits a moduli space. For example, there is no scheme $X$ so that $\text{Hom}(\spec k, X)$ consists of the isomorphism classes of $k$-schemes of finite presentation \cite[pp. 207]{gw}. In general, this problem arises from the existence of non-trivial automorphisms of schemes. In a 1959 letter to Serre, Grothendiek remarked ``\ldots the only barrier to utilizing the theory of descent [to represent a moduli problem] is the existence of non-trivial automorphisms.'' This is a case where the moral problems compliment technical ones. Considering schemes up to isomorphism, and hence discarding information about non-trivial automorphisms, should make us uncomfortable. The existence (or lack thereof) of these trivial automorphisms give us crucial information about the symmetries of the geometric objects in question. 
		
		Since non-trivial automorphisms of a scheme $X$ are simultaneously desirable (in that they are a display of symmetry), and cumbersome as described above, we would like to replace our moduli problem $F$ by one that keeps track of these non-trivial automorphisms. This is precisely the motivating idea behind a \emph{stack}. The solution suggests itself, at least heuristically: we should replace $F\colon\mathbf{Sch}^\text{op} \to \mathbf{Set} $ by a groupoid valued functor $\tilde F\colon\mathbf{Sch}^\text{op} \to \mathbf{Gpd} $. In more technical terms, we are considering, for every scheme $S,$ the fibered category which assigns to $S$ the groupoid of schemes over $S.$ A more precise, but perhaps slightly less intuitive definition, is that a \emph{stack} is a presheaf of groupoids on $\mathbf{Sch}/S$ that satisfies a certain condition called \emph{fpqc descent.} This can be thought of analogously to the ordinary sheaf condition, although it is more technically involved. We will not go into any detail about fpqc descent. We refer the reader to the stacks project \cite{stacks} for more details.  
				
		\end{subsubsection}
		
			\begin{subsection}{Model categories}
			Model categories and $\infty$-categories are structures in which a generalization of the notion of homotopy makes sense. The basic idea of a model category is as a means to effectively solve the \emph{localization problem}; given a complete and cocomplete category $\mathcal{C}$ (i.e. $\mathcal{C}$ contains all small limits and colimits) together with a certain family of morphisms $W$ in $\mathcal{C},$ we want to understand the category $W^{-1}\mathcal{C},$ in which we formally invert the maps in $W.$ Traditionally, these maps are called \emph{weak equivalences}. This terminology is borrowed from the category $\mathbf{Top}_*$ of (based) topological spaces, with $W$ being precisely the weak homotopy equivalences. While $W^{-1}\mathcal{C}$ certainly exists, it is in general very difficult to work with in practice. The presence of a $\emph{model structure}$ on $\mathcal{C}$ makes $W^{-1}\mathcal{C}$ easier to analyze and do computations with. With these motivational principles in mind, we present the following definition: 
			\begin{defn}[Model category {\cite[pp. 21-22]{toen14}}]
			A model category consists of the following data:
				\begin{itemize}
				\item A complete and cocomplete category $\mathcal{C};$
				
				\item Three classes -- $W,$ $Fib,$ and $Cof$ -- of maps in $\mathcal{C}$ called \emph{weak equivalences}, \emph{fibrations}, and \emph{cofibrations} respectively. Elements of $W \cap Fib$ (resp.~ $W \cap Cof$) are called \emph{trivial fibrations} (resp.~ \emph{trivial cofibrations}). 
				\end{itemize}
				
			This data is subject to the following axioms.
				\begin{itemize}
					\item{(\emph{Two of three property})} For a triangle
					\[\begin{tikzcd}
					X \ar[rr, "f"] \ar[dr, "gf" swap]&&Y \ar[dl, "g"] \\
					&Z& 
					\end{tikzcd} \in \mathcal{C}\]
					all three of the depicted maps are in $W$ if and only if two of them are. 
				
					\item Maps in $Fib, Cof,$ and $W$ are all stable under composition and retract. 
					
					\item Consider the square
					\[\begin{tikzcd}[column sep = large, row sep = large]
					A \ar[r, "f"] \ar[d, "i" swap]& X \ar[d, "p"]\\
					B \ar[r, "g" swap ] \ar[ur, "h", dashed]& Y
					\end{tikzcd} \in \mathcal{C}\]
					with $i \in Cof$ and $p \in Fib.$ If either $i \in W$ or $p \in W$ then the depicted lift $h$ exists so that both triangles strictly commute. Note that this $h$ is not necessarily unique. 
					
					\item Any morphism $f \colon X \to Y \in \mathcal{C}$ can be factored in two ways: as a trivial cofibration followed by a fibration, and as a cofibration followed by a trivial fibration. These factorizations are required to be functorial in $f.$ 
					
				\end{itemize}
			\end{defn}
			We remark that the completeness and cocompleteness guarantee the existence of an initial object $\emptyset$ and terminal object $*$ in $\mathcal{C}.$ Objects $X \in \mathcal{C}$ so that $\emptyset \xrightarrow{!} X$ is a cofibration (resp. ~$X \xrightarrow{!} *$ is a fibration) are called \emph{cofibrant objects} (resp. ~\emph{fibrant objects}). Furthermore, notice that $W$ contains the isomorphisms, and in particular, the identity morphisms. 
			
			Given a model category structure on $\mathcal{C},$ we wish to study the \emph{homotopy category} $W^{-1}\mathcal{C}$ (also denoted $Ho(\mathcal{C}$)) of $\mathcal{C}.$ We'll need to recall some definitions. Firstly, for $X \in \mathcal{C},$ a \emph{cylinder object} $C(X)$ for $X$ is characterized by factorizing the codiagonal map $X \coprod X \to X$ through a weak equivalence, i.e. 
			\[\begin{tikzcd} X \coprod X \ar[r, "{(i_1, i_0)}"] & C(X) \ar[r, "p", "\sim" swap]& X\end{tikzcd}.\]
			Moreover, $C(X)$ is called a \emph{good cylinder object} for $X$ if the map $(i_1, i_0)$ is a cofibration. We are ready now to define homotopic maps in $\mathcal{C}.$
			
			\begin{defn}
			We say that two maps \begin{tikzcd} X\ar[r, "f", shift left] \ar[r, "g" swap, shift right] & Y\end{tikzcd} are \emph{homotopic} if there exists a diagram
			\[\begin{tikzcd}[row sep = large, column sep = large]
			X \ar[d, "i_1" swap ] \ar[dr, "f"]\\
			C(X) \ar[r, "h"]& Y\\
			X\ar[u, "i_0"] \ar[ur, "g" swap ]
			\end{tikzcd}\]
			such that $C(X)$ is a good cylinder object for $X.$ 
			\end{defn}
		
		A model category structure gives this definition of homotopy significant weight. Notice that $p$ becomes an isomorphism in $Ho(\mathcal{C})$ by definition. By the 2-of-3 property, this implies that $i_1$ and $i_0$ also get sent to isomorphisms; in fact their images in $Ho(\mathcal{C})$ are equal given that they share a common inverse. It follows immediately that $f$ and $g$ have the same imagine in $Ho(C).$ We obtain a functor $\mathcal{C}^{cf}_{/\sim} \to Ho(\mathcal{C})$ from the subcategory $\mathcal{C}^{cf}$ of fibrant and cofibrant objects in $\mathcal{C}$ modulo homotopy equivalence to the homotopy category of $\mathcal{C}.$ We conclude with the following theorem of Quillen, 
		\begin{thm}[{Quillen \cite[\S 1, Thm. 1]{quillen}}]
			The aforementioned functor
			\[\mathcal{C}^{cf}_{/\sim} \to Ho(\mathcal{C})\]
			induces an equivalence of categories. 
		
		\end{thm}
		\end{subsection}
	
	This is a considerably nicer situation. We now can obtain an explicit description of $Ho(\mathcal{C})$ in the presence of a model structure. We should add a disclaimer that, in practice, working within the homotopy category can be cumbersome. In particular, compositions do not lift well from $Ho(\mathcal{C})$ to $\mathcal{C}$ even in the presence of a model structure on $\mathcal{C}.$ This follows from the fact that functors that compose in $Ho(\mathcal{C})$ may compose only up to higher homotopy. In practice, keeping track of these homotopies makes it hopelessly difficult to lift maps from $Ho(\mathcal{C})$ to $\mathcal{C}.$ 
	
	We end by giving the idea of an $\infty-$category, following the notation in \textit{Elements of $\infty$-category theory} by Emily Riehl and Dominic Verity. The basic idea of an $\infty$-category $\mathcal{A}$ is that is keeps track of the \textit{homotopy category} and higher homotopies of maps in $\mathcal{A}.$ We end by recalling the definition of the homotopy category:
	\begin{defn}[{Homotopy category \cite{riehl}}]
		Given an $\infty-$category $\mathcal{A}$ in an $\infty$-cosmos $\mathcal{K},$ the \emph{homotopy category} of $\mathcal{A},$ is the category $h\text{Fun}\left(1, \mathcal{A}\right)$ of homotopy classes of maps $1 \to \mathcal{A}$ in $\mathcal{K}.$
	\end{defn}
	
		\end{subsection}

	\end{section}
	
	\begin{section}{Motivation: Serre's Intersection Theorem}
	Derived algebraic geometry becomes useful when considering particularly ``bad'' geometric objects. For example, say we are given irreducible algebraic (i.e. ~generated by systems of polynomial equations) subsets $Y$ and $Z$ inside some smooth algebraic variety $X$ (i.e. ~a sufficiently nice scheme). Given some irreducible component $W$ of $Y\cap Z,$ we'd like to know the number $i(X,Y.Z, W)$ of times that $Y$ and $Z$ meet along $W.$  This flavor of question is the starting point of intersection theory, a well-studied branch of algebraic geometry. Now, if $Y$ and $Z$ meet transversally (i.e. ~the tangent spaces of $Y$ and $Z$ generate all of $X$) then the answer is classically understood to be $1.$ On the other hand, if the intersection of $Y$ and $Z$ is not transversal, we run into problems. Heuristically, we can think of a non-transverse intersection as some sort of ``tangling'' of $Y$ and $Z$ within $X.$ In this case, the problem is resolved by a remarkable result, which we present below without proof.
	
	\begin{thm}[{Serre's Intersection Formula \cite[pp. 7]{toen14}}]
		Given schemes $X, Y, Z, W$ as in the above discussion. Under the added assumption of tor independence, the number of intersections of $Y$ and $Z$ in $W$ is given by the expression 
		\[i(X,Y.Z,W) = \text{length}_{\mathscr{O}_{X, W}}\left(\mathscr{O}_{Y, W} \otimes_{\mathscr{O}_{X, W}} \mathscr{O}_{Z, W}\right) \tag{\dag}.\]
		In the general case, we have
		\[i(X,Y.Z,W) = \sum_{n \in \N}(-1)^n\text{length}_{\mathscr{O}_{X, W}}\left(\text{Tor}_n^{\mathscr{O}_{X,W}}\left(\mathscr{O}_{Y, W}, \mathscr{O}_{Z, W}\right)\right).\]
		
	\end{thm}
		
	There is much to unpack in the above formula. A key fact is that $(\dag)$ may be obtained directly from the schematic intersection $Y \cap Z.$ That is, if $Y$ and $Z$ intersect in a sufficiently nice way along $W,$ then the number of intersections can be ascertained solely from the scheme $Y \cap Z.$ On the other hand, in the case that $Y$ and $Z$ meet pathologically, information about the intersection is hidden in higher degree homological data \cite{toen14}. In this case, the role of derived algebraic geometry is to replace the scheme $Y \cap Z$ with a \emph{derived scheme} that will naturally contain the higher-order information originally hidden within the tor-algebras.

	\end{section}

	\begin{section}{Derived algebraic geometry}
		This section largely follows the method proposed in section two of \cite{toen14}. Information on the theory of simplicial rings was obtained from various sources, primarily \cite{akhil} and \cite{quillen}. For further reading on derived stacks see \cite{vez06}. For a different development of derived schemes using $E_\infty$ spaces, see \cite{lurie04}. 
		\begin{subsection}{Simplicial rings}
		 
		The above discussions suggest that we aught to augment a scheme by allowing for a notion of homotopy. As schemes are constructed from commutative rings, it is natural to first find a suitable replacement for $\mathbf{Comm},$ and then build a theory of derived schemes from this replacement. A natural choice is to replace a ring with a simplicial ring, i.e., a simplicial object in the category of rings. To be precise, a simplicial ring $A_\bullet$ consists of a sequence $\{A_i\}_{i = 0}^\infty$ of rings together with face and degeneracy ring homomorphisms as depicted below
		\[	\begin{tikzcd}
				 A_0 	\ar[r, "s_0" ]& A_1\ar[l, shift left = 2ex, "d_0" swap]\ar[l, shift right = 2ex, "d_1"swap]  \ar[r, shift right = 1ex] \ar[r, shift left = 1ex]&\ar[l, shift left = 2ex, ] \ar[l] \ar[l, shift right = 2ex]  \cdots.
			\end{tikzcd}
		\]	
		Morphisms of simplicial rings are morphisms of simplicial sets so that the levelwise maps are ring homomorphisms. In this way we obtain a strict $1$-category $\mathbf{sComm}$ of simplicial commutative rings.  
		\begin{eg}[{Free simplicial rings \cite[pp. 2]{akhil}}]
			Let $X_\bullet \in \mathbf{sSet}$ and $R \in \mathbf{Comm}.$ Then we have $R[X_\bullet]$ is a simplicial commutative ring, where the $n^\text{th}$ degree is given by the polynomial ring $R[X_n].$ Face and degeneracy homomorphisms on $R[X_\bullet]$  are inherited from the face and degeneracy maps on the $X_i,$ e.g. if $d\colon X_n \to X_{n-1}$ is a face map, then $\tilde d \colon R[X_n] \to R[X_{n-1}]$ is the homomorphism obtained by $x \in X_n \mapsto d(x) \in X_{n-1}$ and $1 \mapsto 1.$ 
		\end{eg}
		
		By a 1954 theorem of Moore{\cite[\S 3 Thm. 3]{moore55}}, we know that the underlying simplicial set of a simplicial group is a Kan complex. By the existence of forgetful functors $\mathbf{sComm} \to \mathbf{sGrp} \to \mathbf{sSet},$ it follows that the underlying simplicial set of $A_\bullet \in \mathbf{sComm}$ is also a Kan complex. This allows us to define homotopy groups $\pi_*A_\bullet$ of a simplicial ring $A_\bullet$ (the basepoint is taken to be $0$) as the homotopy groups of the underlying Kan complexes. It turns out that $\pi_*$ is a functor from $\mathbf{sComm}$ to the category $\mathbf{GrComm_{\geq 0}}$ of (nonnegatively) graded commutative rings.\cite[\S 1 Proposition 1]{akhil} 
		
		Recall from section $1.1.1.$ that Grothendiek's notion of a scheme was a generalization of classical algebraic geometry (which studied finitely generated, reduced algebras over an algebraically closed field) to arbitrary commutative rings. In particular, the notion of a scheme allows us to study arbitrary rings, including those that have nilpotents, geometrically. Derived algebraic geometry continues this story by allowing for ``extra nilpotents'' in higher homotopy groups. Therefore, we would expect there to be an analogous relation in $\mathbf{sComm}$ to that of a ring $R$ with its associated reduced ring $R^\text{red} = R_{/\sqrt{0}}.$ 
		
		To realize this analogy, recall that any ring homomorphism $R \to A$ from a ring $R$ to a reduced ring $A$ factors uniquely through $R^\text{red}.$ In other words, the category of reduced rings is a reflective subcategory $\mathbf{Comm}$ . Notice that there is a fully faithful embedding $\mathbf{Comm} \hookrightarrow \mathbf{sComm}$ given by sending a ring $R$ to the constant simplicial ring $R.$ This inclusion admits a left adjoint, namely $\pi_0 \colon \mathbf{sComm} \to \mathbf{Comm}.$ Hence we can view a simplicial ring $A_\bullet$ as an augmentation of $\pi_0(A_\bullet)$ that allows for the addition of higher dimensional homotopical data. Moreover, as suggested by the discussion in section 2.1, the higher homotopy rings of $A_\bullet$ \emph{should} be thought of as a subtler notion of nilpotents. This notion will be made more precise in the next section, but a suggestive fact in this direction is that all simplicial fields are constant.
		
		Our goal was to develop an extension of $\mathbf{Comm}$ with homotopical properties. While we have defined homotopy groups of objects in $\mathbf{sComm},$ the picture is incomplete in that the above discussion has been 1-categorical in nature, and thus does not allow for a notion homotopy between maps $A_* \to B_* \in \mathbf{sComm}.$  The solution is to simplicially enrich $\mathbf{sComm}.$ For some for some $X \in \mathbf{sSet}$ and $A_\bullet \in \mathbf{sComm}$ we define the tensor $X_\bullet \otimes A_\bullet$ as  
		\[\left(X_\bullet \otimes A_\bullet\right)_n = \bigotimes_{X_n} A_n,\]
		where we are taking the tensor product of $A_n$ with itself indexed by elements of $X_n.$ 
		Now for $A_\bullet, B_\bullet \in \mathbf{sComm},$ we define $\Fun(A_\bullet, B_\bullet) \in \mathbf{sSet}$ as
		\[\left(\Fun(A_\bullet, B_\bullet)\right)_n = \text{Hom}_\mathbf{sComm}\left(\Delta^n \otimes A_\bullet, B_\bullet\right).\]
		We note that this process is completely analogous to the simplicial enrichment of $\mathbf{sSet}$ and holds if we replace $\mathbf{sComm}$ with $\mathcal{C}^{\mathbf{\Delta}^\text{op}}$ for any cocomplete category $\mathcal{C}.$
		
		We have now realized $\mathbf{sComm}$ as an $\infty$-category, i.e. a category enriched over $\mathbf{sSet}.$ This allows us to easily define homotopic maps! Indeed, we say that two maps
		\[\begin{tikzcd} A_\bullet \ar[r, shift left, "f"] \ar[r, shift right, "g"swap]& B_\bullet \end{tikzcd} \in \Fun\left(A_\bullet, B_\bullet\right)\]
		are (simplicially) $\emph{homotopic}$ if there exists a map 
		\[\begin{tikzcd} \Delta^1 \otimes A_\bullet \ar[r] & B_\bullet\end{tikzcd}\]
		whose restriction to $\Delta^{\{0\}} \otimes A_\bullet$ (resp. $\Delta^{\{1\}} \otimes A_\bullet$) is $f$ (resp. $g$); compare this with definition 1.3. 
		
		Before concluding this section, a brief digression on $\infty$-categories is in order. We have implicitly used the term $\infty$-category above to mean a simplicially enriched category. This is a somewhat stricter notion than is standard in the sense that composition along $0$-cells associative in unital in the strong sense (i.e., $0$-cells compose strictly, not up to some higher homotopy). In the literature, $\infty$-categories often refer to quasi-categories (i.e. simplicial sets with the inner horn lifting property). If $\mathbf{sComm}$ were already enriched over Kan complexes, then we could apply the homotopy coherent nerve construction to obtain an associated quasi-category $N\left(\mathbf{sComm}\right).$ Unfortunately, this is not the case. On the other hand, it is well known that a simplicial model structure gives a model for $(\infty, 1)$-categories. A theorem of Quillen gives us precisely what we need:
		\begin{thm}[{Quillen \cite{quillen}}]
		There exists a cofibrantly generated simplicial model structure on $\mathbf{sComm}$ with weak equivalences and fibrations being those maps that induce weak equivalences and fibrations on the underlying simplicial sets. One obtains a Quillen adjunction between $\mathbf{sSet}$ and $\mathbf{sComm}$ induced by the free-forgetful adjunction on the underlying 1-categories.  
		\end{thm}
		
		This theorem allows us to work in the given simplicial model structure. If we prefer to work in quasi-categories, we may consider the subcategory $\mathbf{sComm}^\text{cf}$ of fibrant and cofibrant objects in $\mathbf{sComm},$ which is naturally enriched over Kan complexes. We may then apply the homotopy coherent nerve, giving us an associated quasi-category with an equivalent homotopy category to that of $\mathbf{sComm}.$ 
		\end{subsection}

		\begin{subsection}{Derived schemes}
		We'll begin with the following definition. Recall that a \emph{stack} on an ordinary commutative ring $R$ was a fpqc sheaf of groupoids on $\mathbf{Comm}.$
		
			\begin{defn}
			 Given a topological space $X,$ let $\mathbf{Ouv}_X$ denote the category of open subsets of $X$ with morphisms the inclusions. A \emph{stack of simplicial commutative rings} $\mathscr{O}_X$ on $X$ is a $\mathbf{sComm}$ valued presheaf on $\mathbf{Ouv}_X$ that satisfies fpqc descent. A pair $(X, \mathscr{O}_X),$ where $X$ is a topological space and $\mathscr{O}_X$ is stack of simplicial commutative rings on $X,$ is called a \emph{derived ringed space.}
			\end{defn}
			
			We're now finally ready to define a derived scheme! 
			
			\begin{defn}
			The category $\mathbf{dSch}$ is the full subcategory of derived ringed spaces so that the following two conditions are satisfied:
				\begin{enumerate}
				\item The truncation $(X, \pi_0\left(\mathscr{O}_X\right))$ is a scheme.
				
				\item For all $i,$ the sheaf $\pi_i\left(\mathscr{O}_X\right)$ of $\pi_0\left(\mathscr{O}_X\right)$-modules is quasi-coherent. 
				\end{enumerate}
				
			An \emph{affine derived scheme} is a derived scheme $(X, \mathscr{O}_X)$ so that $(X, \pi_0(X))$ is an affine scheme. The $\infty$-category of derived affine scheme is denoted $\mathbf{dAff}.$ 
			
			\end{defn}
			\begin{prop}[To\"en {\cite[pp.32]{toen14}}]
			There is an equivalence of $\infty$-categories $\mathbf{sComm}^\text{op} \to \mathbf{dAff}.$
			\end{prop}
			Recall that we have an adjunction $\pi_0 \dashv i$ of functors 
			\[\begin{tikzcd} \mathbf{sComm} \ar[r, shift left, "\pi_0"]& \mathbf{Comm} \ar[l, shift left, "i"].\end{tikzcd}\]
			This extends to an adjunction $t_0 \dashv i,$ on the level of derived schemes
			\[\begin{tikzcd} \mathbf{dSch} \ar[r, shift left, "t_0"]& \mathbf{Sch} \ar[l, shift left, "i"].\end{tikzcd}\]
			where $t_0(X, \mathscr{O}_X) = (X, \pi_0\left(\mathscr{O}_X\right))$ is the truncation. We see then that $\mathbf{Sch}$ sits in $\mathbf{dSch}$ as a reflective subcategory. 
			
			For every $X \in \mathbf{dSch}$ the counit \begin{tikzcd}[column sep = small]\epsilon\colon t_0\circ i \ar[r, Rightarrow] &\text{id}_\mathbf{dSch}\end{tikzcd} of the adjunction gives a map
			\[\begin{tikzcd}j \colon t_0(X) \ar[r] & X \end{tikzcd}\]
			 that is natural in $X.$ This should be thought of as the analogue of the inclusion
			 \[\begin{tikzcd}Y_{red} \ar[r, hook]& Y\end{tikzcd} \in \mathbf{Sch}\]
			of the underlying reduced scheme $Y_{red}$ of a scheme $Y.${\cite[pp.31]{toen14}} Therefore, the derived scheme $X$ can be thought of as an infinitesimal thickening of $t_0(X)$ via higher dimensional homotopic data. The truncation $t_0$ has a generalization. We have a Postnikov tower
			\[\begin{tikzcd}
			\mathscr{O}_X \ar[d] \\ 
			\vdots \ar[d] \\ 
			t_n \left(\mathscr{O}_X\right)\ar[d] \\
			t_{\leq n-1}\left(\mathscr{O}_X\right) \ar[d]\\ 
			\vdots \ar[d] \\
			t_0\left(\mathscr{O}_X\right) = \pi_0\left(\mathscr{O}_X\right)
			\end{tikzcd} \tag{$*$}\] 
			of derived schemes characterized by the following two properties\cite[pp.33]{toen14}:
			\begin{enumerate}
			\item For all $i \geq n, \pi_i\left(t_{\leq n}\left(\mathscr{O}_X\right)\right) = 0.$
			
			\item For all $i \leq n,$ the morphism $\mathscr{O}_X \rightarrow t_n\left(\mathscr{O}_X\right)$ induces isomorphisms \[\pi_i\left(\mathscr{O}_X\right)\cong \pi_i\left(t_{\leq n}\left(\mathscr{O}_X\right)\right).\] 
			\end{enumerate}
			For a derived scheme $(X, \mathscr{O}_X)$ we denote by $t_{\leq n}(X)$ the derived scheme $(X, t_{\leq n}\left(\mathscr{O}_X\right)).$ Then the tower $(*)$ induces a diagram of derived schemes 
			\[\begin{tikzcd}t_0(X) \ar[r] & \cdots \ar[r] & t_{\leq n}(X) \ar[r] &\cdots \ar[r] &X \end{tikzcd}\]
			so that $X = \text{hocolim}_n t_{\leq n}(X).$ In this way, we may interpolate between schemes and derived schemes by a sequence $\mathbf{dSch}_{\leq n}$ of reflective subcategories of $\mathbf{dSch},$ with reflectors given by the $t_{\leq n}.$ This allows for simplified descriptions of mapping spaces between derived schemes\cite{toen14}:
			\[Map(X, Y) \simeq \text{holim}_nMap(t_{\leq n}\left(X\right),t_{\leq n}\left(Y\right)).\]

			Throughout this paper, we have proposed a point of view of derived schemes as a means to package higher dimensional homotopitcal ``nilpotents'' of schemes. Giving a complete justification is beyond the scope of the current project, but we end by giving an example that points in this direction. 
			
			\begin{eg}[{Fiber products \cite[pp.33]{toen14}}]
			By the equivalence of $\infty$-categories \\$\mathbf{sComm}^\text{op} \to \mathbf{dAff},$ a diagram  
			\[\begin{tikzcd} X \ar[r]& Z  & Y\ar[l]\end{tikzcd}\]
			of affine derived schemes corresponds to a diagram 
			\[\begin{tikzcd} A & C\ar[r]\ar[l] & B\end{tikzcd}\]
			of simplicial commutative rings. In the case that $A, B,$ and $C$ are rings (regarded as constant simplicial rings), the (homotopy) pushout of the above diagram is a simplical ring $D$ so that $\pi_n(D) \cong Tor_n^C(A, B).$ In general, the fiber product $X \times_Z Y$ of derived schemes can be computed by gluing together the affine pieces described above. In the case that $X, Y,$ and $Z$ are schemes, $X \times_Z Y$ is a derived scheme whose truncation $t_0\left(X \times_Z Y\right)$ is the ordinary fiber product of schemes. Moreover, $\pi_n\left(X \times_Z Y\right) \cong Tor_n^{\mathscr{O}_Z}\left(\mathscr{O}_X, \mathscr{O}_Y\right).$  We see now a connection with the Serre intersection formula presented in section 2, and how higher tor-algebras fall naturally out of the theory of derived schemes. 
			\end{eg}
		
		\end{subsection}

	\end{section}

\newpage

\bibliography{DAG}

\end{document}